\documentclass{amsart}
\usepackage{amsfonts,amssymb,amscd,amsmath,enumerate,verbatim}
\usepackage[all]{xy}
\SelectTips{eu}{}
\xyoption{curve}

\newcommand{\xra}{\xrightarrow}

\newcommand{\card}{{\operatorname{card}}}

\newcommand{\wh}{\widehat}
\newcommand{\wt}{\widetilde}

\newcommand{\les}{{\scriptscriptstyle\leqslant}}
\newcommand{\ges}{{\scriptscriptstyle\geqslant}}

\newcommand{\vgu}{{\bigcup\!\scriptstyle{\mathbf\vert}}}

\newcommand{\dcat}[1]{{\mathsf{D}(#1)}}
\newcommand{\dcatf}[1]{{\mathsf{D}^{\mathsf f}(#1)}}
\newcommand{\lotimes}[1]{\otimes^{\bf L}_{#1}}

\newcommand{\Rhom}[3]{{\mathsf R}\!\Hom{#1}{#2}{#3}}

\newcommand{\shift}{{\sf\Sigma}}
\newcommand\col{\colon}

\newcommand\dd{\partial}
\newcommand{\hh}[1]{\operatorname{H}(#1)}
\newcommand{\HH}[2]{\operatorname{H}_{#1}(#2)}
\newcommand{\CH}[2]{\operatorname{H}^{#1}(#2)}

\newcommand\ZZ{\operatorname{Z}}

\newcommand{\Tor}[4]{\operatorname{Tor}_{#1}^{#2}(#3,#4)}
\newcommand{\Ext}[4]{\operatorname{Ext}^{#1}_{#2}(#3,#4)}
\newcommand{\Hom}[3]{\operatorname{Hom}_{#1}(#2,#3)}

\newcommand\Ker{\operatorname{Ker}}
\newcommand{\Ann}{\operatorname{Ann}}

\newcommand\codim{\operatorname{codim}}

\newcommand\edim{\operatorname{edim}}

\newcommand\height{\operatorname{height}}

\newcommand\length{\operatorname{length}}
\newcommand{\lol}[2]{{\operatorname{\ell\hspace{-1pt}\ell}}_{#1}#2}
\newcommand{\lolk}[2]{{\operatorname{\ell\hspace{-1pt}\ell}}^{k}_{#1}#2}

\newcommand{\rank}{\operatorname{rank}}

\newcommand\idmap{\operatorname{id}}

\newcommand\pd{\operatorname{proj\,dim}}

\newcommand\gldim{\operatorname{gl\,dim}}

\newcommand\fm{{\mathfrak m}}
\newcommand\fn{{\mathfrak n}}

\newcommand\fq{{\mathfrak q}}

\newcommand\bsq{{\boldsymbol q}}

\newcommand\BN{{\mathbb N}}

\newcommand\BR{{\mathbb R}}
\newcommand\BZ{{\mathbb Z}}

\newcommand\ca{{\mathsf A}}

\newcommand\cc{{\mathsf C}}

\newcommand\SH{{\mathsf h}}

\newcommand\SJ{{\mathsf j}}
\newcommand\SK{{\mathsf k}}

\newcommand\ST{{\mathsf t}}
\newcommand{\eps}{{\varepsilon}}
\newcommand{\vf}{{\varphi}}

\newcommand{\var}{{\hskip1pt\vert\hskip1pt}}

\theoremstyle{plain}
\newtheorem{theorem}{Theorem}[section]

{\Alph{theorem}}
\newtheorem{proposition}[theorem]{Proposition}
\newtheorem{lemma}[theorem]{Lemma}

\theoremstyle{definition}

\newtheorem{construction}[theorem]{Construction}
\newtheorem{example}[theorem]{Example}

\newtheorem{chunk}[theorem]{}

\theoremstyle{remark}

\newtheorem{remark}[theorem]{Remark}
\numberwithin{equation}{theorem}

\newcommand{\cxy}[2][R]{{\operatorname{cx}_{#1}{#2}}}
\newcommand{\injcxy}[2][R]{{\operatorname{inj\,cx}_{#1}{#2}}}
\newcommand{\level}[3]{\operatorname{level}_{#1}^{#2}(#3)}
\newcommand{\thickn}[3]{{\mathsf{thick}}_{#2}^{#1}(#3)}
\newcommand{\thick}[2]{{\mathsf{thick}}_{#1}(#2)}

\begin{document}
\title[Cohomology over complete intersections]
{Cohomology over complete intersections\\
via exterior algebras}

\date{\today}

\author[L.~L.~Avramov]{Luchezar L.~Avramov}
\address{Department of Mathematics, University of Nebraska, Lincoln, NE 68588, U.S.A.}
\email{avramov@math.unl.edu}
\author[S.~B.~Iyengar]{Srikanth B.~Iyengar} 
\address{Department of Mathematics, University of Nebraska, Lincoln, NE 68588, U.S.A.}  
\email{iyengar@math.unl.edu}

\thanks {Research partly supported by NSF grants DMS 0803082 (L.L.A)
and DMS 0602498 (S.B.I).  Part of the work of the second author was done
at the University of Paderborn on a visit supported by a For\-schungs\-preis
from the Humboldt\-stiftung.}

  \begin{abstract}
A general method for establishing results over a commutative complete 
intersection local ring by passing to differential graded modules over 
a graded exterior algebra is described.  It is used to deduce, in a uniform 
way, results on the growth of resolutions 
of complexes over such local rings. 
  \end{abstract}

\dedicatory{To Karin Erdmann on her 60th birthday.}


\maketitle

\section*{Introduction}

This paper concerns homological invariants of modules and complexes over 
complete intersection local rings. The goal is to explain a method by which one 
can establish in a uniform way results over such rings by deducing them from
results on DG (that is, differential graded) modules over a graded exterior algebra,
which are often easier to prove. A secondary purpose is to demonstrate
the use of numerical invariants of objects in derived categories, called 
`levels', introduced in earlier joint work with Buchweitz and Miller~\cite{ABIM};
see Section~\ref{Levels}.  Levels allow one to 
track homological \emph{and} structural information under changes of rings 
or DG algebras, such as those involved when passing from complete 
intersections to exterior algebras.

We focus on the complexity and the injective complexity of a complex 
$M$ over a complete intersection ring $R$. These numbers measure, on a 
polynomial scale, the rate of growth of the minimal free resolution and the 
minimal injective resolution of $M$, respectively.  The relevant basic 
properties are established in Section~\ref{Complexities}.

Complexities can be expressed as dimensions of certain algebraic varieties, 
attached to $M$ in \cite{Av:vpd}, and earlier proofs of key results relied on 
that theory.  In Section~\ref{Complete intersection local rings.I} we deduce 
them from results on differential graded Hopf algebras, presented with 
complete proofs in Section~\ref{Differential Graded Hopf algebras}.  These
pleasingly simple proofs build on nothing more than basic homological 
algebra, summarized in Section~\ref{Composition products}.

In the last three sections of this article the goal is to link the complexity 
of $M$ and the Loewy length of its homology modules.  In~\cite{ABIM} 
such a result is deduced from a more general statement, which applies 
to arbitrary local rings.  From that paper we import DG versions of the 
New Intersection Theorem, recalled in Section~\ref{Projective levels}, 
and of the Bernstein-Gelfand-Gelfand correspondence, which we
refine in Section~\ref{Exterior algebras.II}.  This allows us to establish 
a result on complexities over exterior algebras, which we then translate 
in Section~\ref{Complete intersection local rings.II} into the desired link 
between Loewy length and complexity over complete intersection local rings.
 
\section{Levels}
\label{Levels}

In this section $A$ denotes a DG algebra.  

We write $\dcat A$ for the derived category of DG $A$-modules;
it is a triangulated category with shift functor denoted $\shift$.
The underlying graded object of a DG object $X$ is denoted $X^\natural$.
Rings are treated as DG algebras concentrated in degree zero; over a ring
DG modules are simply complexes, and modules are DG modules concentrated
in degree zero; see \cite[\S3]{ABIM} for more details and references.  Unless
specified otherwise, all DG modules have left actions.

A non-empty subcategory $\cc$ of $\dcat A$ is said to be \emph{thick} 
if it is an additive subcategory closed under retracts, and every exact triangle 
in $\dcat A$ with two vertices in $\cc$ has its third vertex in $\cc$;
thick subcategories are triangulated.

Given a DG $A$-module $X$, we write $\thick AX$ for the smallest thick subcategory 
containing $X$.  The existence of such a subcategory can be seen by realizing
it as the intersection of all thick subcategories of $\dcat A$ that contain $X$.  
Alternatively, the objects of $\thick AX$ can be built from $X$ in a series of steps,
described below.

\begin{chunk}
\label{chu:levels}
For every $X$ in $\dcat A$ and each $n\ge0$ we define a full subcategory 
$\thickn n{A}{X}$ of $\thick AX$, called the $n$th \emph{thickening} of $X$ in 
$\ca$, as follows. Set $\thickn 0{A}{X}=\{0\}$; the objects of $\thickn 1{A}{X}$ are 
the retracts of finite direct sums of shifts of $X$. For each $n\ge2$, the objects 
of $\thickn n{A}{X}$ are retracts of those $U\in\dcat A$ that appear in some exact triangle 
$U'\to U\to U''\to\shift U'$ with $U'\in\thickn{n-1}{A}{X}$ and $U''\in\thickn 1{A}{X}$.  

Every thickening of $X$ is clearly embedded in the next one; it is also clear
that their union is a thick subcategory of $\ca$ containing $X$,
which is therefore equal to $\thick A{X}$.  Thus, $\thick A{X}$
is equipped with a filtration
\[
\{0\}=\thickn 0{A}{X}\subseteq \thickn 1{A}{X}\subseteq\cdots
\subseteq\bigcup_{n\in\BN} \thickn n{A}{X} = \thick A{X}\,.
\]

To each object $U$ in $\dcat A$ we associate the number
\[
\level{A}{X}U = \inf\{n\in \BN\mid U\in \thickn n{A}{X}\}
\]
and call it the \emph{$X$-level of $U$ in $\dcat A$}.  It measures
the number of extensions needed to build $U$ out of $X$.  Evidently,
$\level{A}{X}U<\infty $ is equivalent to $U\in\thick A{X}$.
  \end{chunk}

We refer the reader to \cite[\S\S2--6]{ABIM} for a systematic study of
levels.  The  properties used explicitly in this paper are recorded below.

\begin{chunk}
\label{ch:levels}
Let $U$, $X$, and $Y$ be objects in $\dcat A$.
\begin{enumerate}[\rm(1)]
\item
If $\thickn 1{A}X=\thickn 1{A}Y$, then one has 
$\thickn n{A}X=\thickn n{A}Y$ for all $n$, and 
hence an equality $\level{A}XU = \level{A}YU$.
\item
If $B$ is a DG algebra and $\SJ\col\dcat A\to\dcat B$ is an exact 
functor, then
\[
\level{A}XU \geq \level{B}{\SJ(X)}{\SJ(U)}\,.
\]
Equality holds when $\SJ$ is an equivalence.
\end{enumerate}
 \end{chunk}

A level of interest in this paper is related to the notion of Loewy length.

\begin{chunk}
  \label{ch:loewy}
Let $B$ be a ring and $k$ a simple $B$-module. For each $B$-module 
$H$, $\lolk BH$ denotes the smallest  integer $l\ge 0$ such that $H$ has 
a filtration by $B$-submodules
  \[
0=H^0\subseteq H^1 \subseteq\cdots\subseteq H^{l-1}\subseteq H^{l}=H
  \]
with every $H^i/H^{i-1}$ isomorphic to a sum of
copies of $k$; if no such filtration exists, we set $\lolk BH=\infty$.  
If $B$ has a unique maximal left ideal $\fm$, then $k$ is isomorphic
to $B/\fm$ and $\lolk BH$ is equal to the \emph{Loewy length}
of $H$, defined by the formula
\[
\lol BH = \inf\{l\in\BN\mid \fm^{l}H=0\}\,;
\]
see \cite[6.1.3]{ABIM}.  If $\length_{B}H$ is finite, then so is 
$\lol BH$, and the converse holds when $H$ is noetherian; see 
\cite[6.2(4)]{ABIM}.
\end{chunk}

We say that the DG algebra $A$ is \emph{non-negative} if $A_i=0$ for $i<0$.  
When this is the case, there is a canonical morphism of DG algebras 
$A\to\HH0A$, called the \emph{augmentation} of $A$; it turns every 
$\HH0A$-module into a DG $A$-module.  

\begin{chunk}
\label{ch:klevel}
When $A$ is non-negative and $k$ is a simple $\HH0A$-module
one has:
\begin{enumerate}[{\rm(1)}]
\item
The number $\level AkU$ is finite if and only if the $\HH0A$-module
$\bigoplus_{n\in\BZ}\HH nU$ admits a finite filtration with subfactors
isomorphic to $k$.
\item
When $\level AkU$ is finite the following inequality holds:
\[
\level AkU\leq \sum_{n\in\BZ} \lolk{\HH0A}{\HH nU}\,.
\]
\item
When $\HH 0A=k$ and the $k$-module $\hh U$ is finite, one has an inequality 
\[
\level AkU\leq \card\{n \in\BZ\mid \HH nU\ne 0\}\,.
\]
\end{enumerate}
Indeed, (1) and (2) are contained in \cite[6.2(3)]{ABIM}, while (3) is extracted 
from \cite[6.4]{ABIM}.
\end{chunk}

\section{Complexities}
\label{Complexities}

In this section we introduce a notion of complexity for DG modules over 
a suitable class of DG algebras, and establish some of its elementary 
properties.  We begin with a reminder of the construction of derived 
functors on the derived category of DG modules over a DG algebra; 
see \cite[\S1]{AH} for details.

  \begin{chunk}
\label{ch:Ext}
Let $A$ be a DG algebra and $A^{\mathsf o}$ its opposite DG algebra. 

A \emph{semifree filtration} of a DG $A$-module $F$ is a filtration 
  \[
0=F^0\subseteq F^1 \subseteq\cdots\subseteq F^{n-1}\subseteq F^{n}
\subseteq\cdots
  \]
by DG submodules with $\bigcup_{n\ges0}F^{n}=F$ and each DG module
$F^n/F^{n-1}$ isomorphic to a direct sum of suspensions of $A$; when
one exists $F$ is said to be semifree.  When $F$ is semifree its underlying
graded $A^\natural$-module $F^\natural$ is free, and the functors 
$\Hom AF-$ and $(-\otimes_AF)$ preserve quasi-isomorphisms of DG 
$A$-modules.

A \emph{semifree resolution} of a DG module $U$ is a quasi-isomorphism $F\to U$ of DG modules with $F$ semifree; such a resolution always exists. If $U\to V$ is a quasi-isomorphism of DG $A$-modules and $G\to V$ is a semifree resolution, then there is a unique up to homotopy morphism $F\to G$ of DG $A$-modules such that the composed maps $F\to G\to V$ and $F\to U\to V$ are homotopic.

These properties imply that after choosing a resolution $F_U\to  U$ for
each $U$, the assignments $(U,V)\mapsto\Hom A{F_U}V$ and $(W,U)\mapsto
W\otimes_A{F_U}$ define exact functors
  \[
\Rhom A--\col\dcat A^{\mathsf o}\times\dcat A\to\dcat\BZ
  \quad\text{and}\quad
-\lotimes A-\col\dcat{A^{\mathsf o}}\times\dcat A\to\dcat\BZ\,,
  \]
respectively.  In homology, they define graded abelian groups
  \[
\Ext{}AUV=\hh{\Rhom AUV}
  \quad\text{and}\quad
\Tor{}AWU=\hh{W\lotimes A U}\,.
  \]
In case of modules over rings these are the classical objects.
  \end{chunk}

We want to measure, on a polynomial scale, how the `size' of $\Ext nAUV$
grows when $n$ goes to infinity.  In order to do this we use the notion of
\emph{complexity} of a sequence $(b_n)_{n\in\BZ}$ of non-negative
numbers, defined by the following equality
  \[
\cxy[{}]{(b_n)}:=
\inf\left\{d \in\BN \left|
\begin{gathered}
\text{there is a number $a\in\BR$ such that}\\
b_n\leq a n^{d-1} \text{ holds for all $n\gg0$}
\end{gathered}
\right\}\right. .
  \]

Throughout the rest of this section,  $(A,\fm,k)$ denotes a \emph{local DG
algebra}, by which we mean that $A$ is non-negative, $A_0$ is a noetherian
ring contained in the center of $A$, $\fm$ is the unique maximal ideal
of $A_0$, and $k$ is the field $A_{0}/\fm$.

The \emph{complexity} of a pair $(U,V)$ of DG $A$-modules is the number
  \[
\cxy[A]{(U,V)} =
\cxy[{}]{\big(\rank_{k}(\Ext nAUV\otimes_{A_0}k)\big)}\,.
  \]
The \emph{complexity} and the \emph{injective complexity} of $U$
are defined, respectively, by
\begin{align*}
\cxy[A]U= \cxy[A]{(U,k)} 
   \quad\text{and}\quad
\injcxy[A]U= \cxy[A]{(k,U)}\,.
\end{align*}
When $A$ is a local ring $\cxy[A]U$ is the polynomial rate of growth of 
$\rank_AF_n$, where $F$ is a minimal free resolution of $U$, while $\injcxy[A]U$ 
is that of the multiplicity of an injective envelope of $k$ in $I^n$, where $I$
is a minimal injective resolution of $U$.

The following properties of these invariants have been observed before
for modules, and for complexes, over local rings.  We extend them to handle 
DG modules over DG algebras, sometimes with alternative proofs.

\begin{lemma}
\label{lem:cxy1}
If $U$ is in $\thick AX$ for some DG $A$-module $X$, then one has
\[
\cxy[A]U\le \cxy[A]X\quad\text{and}\quad  
\injcxy[A]U\le \injcxy[A]X\,.
\]
Equalities hold in case $\thick {A}U=\thick {A}X$.
\end{lemma}

\begin{proof}
We verify the inequality for complexities; a symmetric argument yields
the one for injective complexities.

The number $h=\level AXU$ is finite; we induce on it.  When $h$ equals one $U$ is a retract of $\bigoplus_{j=1}^s\shift^{i_j}X$ for some integers $i_1,\dots,i_s$, so $\Ext nAUk$ is a direct summand of $\bigoplus_{j=1}^s\Ext{n+i_j}AXk$, and the desired inequality is clear. For $h\ge1$, one may assume there is  an exact triangle $U'\to U\to U''\to\shift U'$ with $\level AX{U'}\le h-1$ and $\level AX{U''}=1$; the associated cohomology
exact sequence shows that the induction hypothesis implies the desired
inequality for complexities.
\end{proof}

\begin{lemma}
\label{lem:cxy2}
If $F$ is a finite free complex of $A_0$-modules with $\hh F\ne 0$, 
then 
\[
\cxy[A]{(U\otimes_{A_0}F)}=\cxy[A]U\quad\text{and}\quad
\injcxy[A]{(U\otimes_{A_0}F)}=\injcxy[A]U\,.
\]
\end{lemma}

\begin{proof}
Since $F$ is a finite free complex, in $\dcat A$ one has isomorphisms
\[
\Rhom Ak{U\otimes_{A_0}F}\simeq \Rhom AkU \lotimes{A_0}F \simeq 
\Rhom AkU \lotimes k(k\otimes_{A_0}F)\,.
\]
This yields, for each integer $n$, an isomorphism of $k$-vectorspaces
\[
\Ext nAk{U\otimes_{A_0}F}\cong \bigoplus_{i\in\ZZ}\Ext {n+i}AkU 
\otimes_{k}\HH {i}{k\otimes_{A_0}F}\,. 
\]
Since $\HH i{k\otimes_{A_{0}}F}$ is finite for each $i$, and zero for $|i|\gg0$ but not for all $i$,
the equality of injective complexities follows. The argument for complexities is similar.
  \end{proof}

\begin{chunk}
\label{chu:quism}
For a morphism of DG algebras $\vf\col A\to B$, let $\vf_*$ be the functor 
forgetting the action of $B$.  One then has an adjoint pair of exact functors
    \[
\xymatrixcolsep{4.5pc}
\xymatrix{
\dcat A\ar@{->}[r]<.6ex>^-{(B\lotimes A-)}
&\dcat B\ar@{->}[l]<.6ex>^-{\vf_*}
}
    \]
that are inverse equivalences if $\vf$ is a quasi-isomorphism; see
\cite[3.3.1, 3.3.2]{ABIM}.
If, furthermore, $U$ is a DG $A$-module, $V$ is a DG $B$-module, 
and $\mu\col U\to V$ is a $\vf$-equivariant quasi-isomorphism, then
by \cite[3.3.3]{ABIM} one has canonical isomorphisms
   \[
U\simeq\vf_*(V)
    \qquad\text{and}\qquad
B\lotimes AU \simeq V\,.
   \]

A DG algebra $B$ is said to be \emph{quasi-isomorphic} to $A$ if there 
exists a chain 
    \[
  \xymatrixcolsep{1.35pc}
  \xymatrix{
A\ar@{->}[r]^-{\simeq}
&A^0\ar@{<-}[r]^-{\simeq}
&A^1\ar@{->}[r]^-{\simeq}
&\ \cdots\ \ar@{<-}[r]^-{\simeq}
&A^{i-1}\ar@{->}[r]^-{\simeq}
&A^{i}\ar@{<-}[r]^-{\simeq}
&B
}
  \]
of quasi-isomorphisms of DG algebras.  Such a chain induces a unique equivalence 
  \[
\SJ\col\dcat A\to\dcat B\,,
  \]
of triangulated categories; if $B_i=0$ for $i<0$, then \cite[3.6]{ABIM} yields
  \[
\SJ(\HH0A)\simeq\HH0B\,.
   \]
    \end{chunk}

\begin{lemma}
\label{lem:quism}
If $(B,\fn,l)$ is a local DG algebra quasi-isomorphic to $A$, then each 
equivalence of derived categories $\SJ\col\dcat A\to\dcat B$ as in 
\emph{\ref{chu:quism}} induces an isomorphism
  \[
\SJ(k)\simeq l
  \]
in $\dcat B$, and for all $U$ and $V$ in $\dcat A$ there is an equality
  \[
\cxy[A]{(U,V)} = \cxy[B]{(\SJ(U),\SJ(V))}\,.
   \]
    \end{lemma}

  \begin{proof}
As $\HH0A$ is equal to $A_0/\dd(A_1)$, it is a local ring with maximal ideal 
$\fm/\dd(A_1)$.  Similarly, $\HH0B=B_0/\dd(B_1)$ is a local ring with maximal 
ideal $\fn/\dd(B_1)$.  The isomorphism of \ref{chu:quism} maps 
$\SJ(\fm/\dd(A_1))$ to $\fn/\dd(B_1)$, and so induces $\SJ(k)\simeq l$.

Since $\SJ$ is an equivalence, in $\dcat B$ there is an isomorphism 
  \[
\SJ(\Rhom AUV)\simeq\Rhom{B}{\SJ U}{\SJ V}\,.
  \]
Passing to homology one obtains an isomorphism of graded modules
  \[
\Ext{}AUV\cong\Ext{}{B}{\SJ U}{\SJ V}
  \]
that is equivariant over the isomorphism $\hh A\cong\hh B$ of graded
algebras induced by the chain of quasi-isomorphisms inducing $\SJ$.  
This yields an isomorphism
  \[
\Ext{}AUV\otimes_{\HH0A}k\cong\Ext{}{B}{\SJ U}{\SJ V}\otimes_{\HH0B}l
  \]
of graded vector spaces that  is equivariant over the isomorphism $k\cong l$.
  \end{proof}

Sometimes there exist an alternative way for computing complexity:

\begin{lemma}
\label{lem:artinian}
If the ring $A_0$ is artinian, then there is an equality 
  \[
\cxy[A]{(U,V)} = \cxy[{}]{(\length_{A_0}\Ext nAUV)}\,.
   \]
    \end{lemma}

  \begin{proof}
Set $r_n=\rank_k(\Ext nAUV\otimes_{A_{0}}k)$. By hypothesis $\fm^s=0$ 
holds for some integer $s$, so by Nakayama's Lemma the $A_0$-module 
$\Ext nAUV$ is minimally  generated by $r_n$ elements.  Thus, one has 
surjective homomorphisms
   \[
k^{r_n}\gets \Ext nAUV\gets A_0^{r_n}
  \]
 of $A_0$-modules.  They yield inequalities
  \[
r_n\le \length_{A_0}\Ext nAUV\le \length_{A_0}(A_0)r_n
  \]
that evidently imply $\cxy[{}]{(r_n)}=\cxy[{}]{(\length_{A_0}\Ext nAUV)}$.
  \end{proof}

\section{Composition products}
\label{Composition products}

Let $A$ denote a DG algebra, and let $U$, $V$, and $W$ be DG $A$-modules.  In this section we recall the construction of products in cohomology of DG $A$-modules.

\begin{construction}
 \label{con:composition}
For all integers $i,j$ there exist \emph{composition products}
  \[
\Ext{j}{A}VW\times\Ext{i}{A}UV\to\Ext{i+j}{A}UW\,.
  \]
Indeed, choose a semifree resolution $F\to U$; thus $\Ext{i}{A}UV=\CH i{\Hom{A}FV}$. An element $[\alpha]\in\Ext{i}{A}UV$ is then the homotopy class of a degree $-i$ chain map $\alpha\col F\to V$ of DG modules.  Similarly, let $[\beta]\in\Ext{j}{A}VW$ be the class of a chain map $\beta\col G\to W$ of degree $-j$, where $G\to V$ is a semifree resolution. As $F$ is semifree and $F\to U$ is a quasi-isomorphism there is 
a unique up to homotopy chain map $\wt\alpha\col F\to G$ whose composition with $G\to V$ is equal to $\alpha$.  One defines the product $[\beta][\alpha]$ in $\Ext{i+j}{A}UW$ to be the class of $\beta\circ\wt\alpha\col F\to W$.

It follows from their construction that composition products are 
associative in the obvious sense.  In particular, $\Ext{}{A}UU$ and 
$\Ext{}{A}VV$ are graded DG algebras, and that $\Ext{}{A}UV$
is a graded left-right $\Ext{}{A}VV$-$\Ext{}{A}UU$-bimodule.
\end{construction}

Let $A$ and $B$ be DG algebras over a commutative ring $k$. Their tensor 
product is a DG algebra, with underlying complex $A\otimes_{k}B$ and 
multiplication defined by
\[
(a\otimes b)\cdot (a'\otimes b')= (-1)^{|b||a'|} (aa'\otimes bb')\,,
\]
where $|a|$ denotes the degree of $a$.  For each DG $B$-module
$X$ the complexes $U\otimes_kX$ and $\Hom kXU$ have canonical 
structures of DG module over $A\otimes_{k}B$, given by
  \begin{align*}
(a\otimes b)\cdot(u\otimes x)&= (-1)^{|b||u|} (au\otimes bx)\,;
\\
\big((a\otimes b)(\gamma)\big)(x)&=(-1)^{|b||\gamma|}a\gamma(bx)\,.
  \end{align*}

\begin{construction}
  \label{con:tensorproducts}
Let $k$ be a field, and let $A$ and $B$ be DG algebras with 
  \[
A_{i}=0=B_{i} \quad\text{for}\quad i<0
\quad\text{and}\quad A_{0}=k=B_{0}\,.
  \]
Let $F\to k$ and $G\to k$ be  semifree resolutions over $A$ and $B$, 
respectively, and  
\[
\omega\col
\Hom AFk\otimes_{k} \Hom BGk\to \Hom{A\otimes_{k}B}{F\otimes_{k}G}k
\]
be the morphism given by $\big(\omega(\alpha\otimes\beta)\big)(a\otimes b)=
(-1)^{|\beta||a|}\alpha(a)\beta(b)$.  Composing $\hh\omega$ with the 
\emph{K\"unneth isomorphism} $[\alpha]\otimes_{k}[\beta]\mapsto
[\alpha\otimes \beta]$ one obtains a homomorphism
  \[
\hh{\Hom AFk}\otimes_{k}\hh{\Hom BGk}\to \hh{\Hom{A\otimes_{k}B}{F\otimes_{k}G}k}
  \]
The map $F\otimes_{k}G\to k\otimes_kk=k$ evidently is a semifree resolution 
over $A\otimes_{k}B$, so the preceding homomorphisms defines a 
\emph{K\"unneth homomorphism}
\[
\kappa\col \Ext{}Akk \otimes_{k}\Ext{}Bkk\to \Ext{}{A\otimes_{k}B}kk \,.
\]
It is follows from the definition that this is a homomorphism of $k$-algebras.
  \end{construction}

\begin{lemma}
\label{lem:Kunneth}
If $\rank_{k}A_{i}$ is finite for each $i$, then $\kappa$ is an isomorphism.
\end{lemma}

\begin{proof}
Under the hypothesis on $A$, one can choose an $F$ that has a semifree filtration 
where  each free graded $A^\natural$-module $(F^{n}/F^{n-1})^\natural$ of finite 
rank; for instance, take $F$ to be the classical bar-construction, see \cite{Mac}.
In this case the map $\omega$ from Construction \ref{con:tensorproducts} is
bijective, so $\hh\omega$ is  an isomorphism.
\end{proof}

\section{Differential Graded Hopf algebras}
\label{Differential Graded Hopf algebras}

In this section $k$ is a field and $A$ is a \emph{DG Hopf algebra} over 
$k$; that is, $A$ is a non-negative DG algebra with $A_0=k$, and with 
a morphism $\Delta\col A\to A\otimes_{k}A$ of DG $k$-algebras, called 
the \emph{comultiplication}, satisfying $(\eps^A\otimes A)\Delta=
\idmap^A=(A\otimes_k\eps^A)\Delta$, where $\eps^{A}\col A\to k$ is the
canonical surjection.

The principal results in this section, Theorems~\ref{thm:dghcx} and 
\ref{thm:symmetry}, establish symmetry properties of complexities. 
The arguments in the proof exploit operations on DG modules 
provided by the Hopf algebra structure on $A$.

  \begin{chunk}
When $A$ is a DG Hopf algebra the DG $(A\otimes_kA)$-modules 
$U\otimes_{k}V$ and $\Hom kUV$ acquire structures of DG 
$A$-modules through $\Delta$.  The isomorphisms
\begin{align}
\label{iso:tensor}
U\otimes_{k}k &\cong U \cong k\otimes_{k}U\\
\label{iso:hom}
\Hom kkV &\cong V\\
\label{iso:adjunction}
\Hom{A}{U\otimes_{k}V}W &\cong \Hom{A}U{\Hom kVW}
\end{align}
are compatible with the respective DG $A$-structures.
  \end{chunk}

\begin{lemma}
\label{lem:retract}
Let $F\to k$ be a semifree resolution and $G$ a semifree DG module. The induced morphism $\pi\col F\otimes_{k}G\to k\otimes_{k}G\cong G$ is a homotopy equivalence.
\end{lemma}

\begin{proof}
The morphism $\pi$ is a quasi-isomorphism because $k$ is a field. As $G$ is semifree, the identity map on $G$ lifts through $\pi$ to give a morphism $\kappa\col G\to F\otimes_{k}G$ such that $\pi\kappa$ is homotopic to $\idmap^{G}$; see \ref{ch:Ext}. Since $F$ is semifree it follows from \eqref{iso:adjunction} that the functor $\Hom{A}{F\otimes_{k}G}-$ preserves quasi-isomorphisms. Noting that $\pi\kappa\pi$ is homotopic to  $\pi\idmap^{F\otimes_{k}G}$, one thus  gets that $\kappa\pi$ is homotopic to $\idmap^{F\otimes_{k}G}$.
\end{proof}

\begin{construction}
Let $A$ be a DG Hopf algebra over $k$ and set $S=\Ext{}Akk$. Let 
$V$ be a DG $A$-module.  Choose semifree resolutions $F\to k$ 
and $G\to V$ over $A$. The assignment $\psi\mapsto \psi\otimes_{k}G$ 
defines a morphism
\[
\Hom{A}FF\to\Hom{A}{F\otimes_{k}G}{F\otimes_{k}G}
  \]
of DG algebras. Since $F\otimes_{k}G$ is homotopy equivalent to $G$, by Lemma~\ref{lem:retract}, in homology it induces a homomorphism of graded $k$-algebras
\[
	\zeta_{V}\col S \to \Ext{}{A}VV\,.
\]
  \end{construction}

The graded center of a graded algebra $B$ consists of the elements $c$
in $B$ that satisfy $bc =(-1)^{|b||c|}cb$ for all $b\in B$.  When every
$c\in B$ has this property $B$ is said to be \emph{graded-commutative}; 
every graded right $B$-module $X$ then has a canonical structure of left
$B$-module, defined by setting $bx=(-1)^{|b||x|}xb$.

\begin{proposition}
\label{prop:centrality}
Let $A$ be a DG Hopf $k$-algebra, and let $U,V$ be DG $A$-modules. 

The algebra $S=\Ext{}Akk$ is then graded-commutative and $\zeta_{V}(S)$
is contained in the graded center of $\Ext{}{A}VV$. Moreover, the
$S$-module structures on $\Ext{}{A}UV$ induced via $\zeta_{V}$ and
$\zeta_{U}$ coincide.
  \end{proposition}

\begin{proof}
We verify the second assertion; the first one follows, as $\zeta_k=\idmap^S$.

We may assume $U$ and $V$ are semifree DG $A$-modules. Let $F\to k$ be a semifree resolution. By Lemma~\ref{lem:retract}, the homology of the complex $\Hom A{F\otimes_{k}U}{F\otimes_{k}V}$ is $\Ext{}AUV$. It implies also that any chain map $F\otimes_{k}U\to F\otimes_{k}V$ is homotopy equivalent to one of the form $F\otimes_{k}\mu$. For any chain map $\sigma\col F\to F$ the compositions
\begin{align*}
&F\otimes_{k}U \xra{\sigma\otimes U} F\otimes_{k} U \xra {F\otimes \mu} F\otimes_{k}V\\
&F\otimes_{k}U \xra{F\otimes \mu} F\otimes_{k} V \xra {\sigma\otimes V} F\otimes_{k}V
\end{align*}
are $\sigma\otimes \mu$ and $(-1)^{|\mu||\sigma|}\sigma\otimes \mu$, respectively. Hence the $S$-actions on $\Ext{}AUV$ induced via $\zeta_{U}$ and $\zeta_{V}$ coincide up to the usual sign. This is the desired result.
\end{proof}

\begin{proposition}
\label{prop:dghfg}
Let $A$ be a DG Hopf $k$-algebra, and $U$ and $V$ be DG $A$-modules.

If the $k$-algebra $S=\Ext{}Akk$ is finitely generated, and the $k$-vector 
spaces $\hh U$ and $\hh V$ are finite, then $\Ext{}{A}UV$ is finite over $S$, 
$\Ext{}{A}UU$, and $\Ext{}{A}VV$.
  \end{proposition}

\begin{proof}
Since $S$ acts on $\Ext{}{A}UV$ through $\Ext{}{A}UU$, finiteness   
over the former implies finiteness over the latter.  The same reasoning
applies with $\Ext{}{A}UU$ replaced by $\Ext{}{A}VV$, so it suffices
to  prove finiteness over $S$.

We claim that the following full subcategory of $\dcat A$ is thick:
\[
\cc = \{W\in\dcat A\var \text{the $S$-module $\Ext{}AWk$ is noetherian}\}\,.
\]
Indeed, it is clear that $\cc$ is closed under retracts and shifts;
furthermore, every exact triangle $W'\to W\to W''\to \shift W'$ in $\dcat A$ 
induces an exact sequence
\[
 \Ext{}A{W''}k \to \Ext{}A{W}k \to \Ext{}A{W'}k
\]
of graded $S$-modules, so if $W'$ and $W''$ are in $\cc$, then so is $W$. 

The condition $\rank_k\hh U<\infty$ implies $U\in\thick Ak$; see 
\ref{ch:klevel}(1).  The finitely generated $k$-algebra $S$ is graded-commutative 
by Proposition~\ref{prop:centrality}, and thus noetherian; this implies $k\in\cc$.  
Now from the thickness of $\cc$ we conclude $U\in\cc$.

{}From a similar argument, now considering the subcategory
\[
\{W\in\dcat A\var \text{the $S$-module $\Ext{}AUW$ is noetherian}\}\,,
\]
one deduces that the $S$-module $\Ext{}AUV$ is finitely generated.
  \end{proof}

The preceding result allows one to draw conclusions about complexities 
from classical results in commutative algebra.  This may not be 
immediately clear to casual users of the subject, as standard textbooks 
leave out important parts of the story by focusing early on on `standard 
graded' algebras.  In the following discussion we refer to Smoke's very 
readable and self-contained exposition in \cite{Sm}.

\begin{chunk}
  \label{ch:HS}
Let $S$ be a graded-commutative $k$-algebra, generated over $k$ 
by finitely many elements of positive degrees, and $C$ a finitely graded
$S$-module.  Replacing generators by their squares, one sees that $C$ 
is also finite over a finitely generated \emph{commutative} $k$-algebra 
$S'$, so there exists a Laurent polynomial $p_C(t)\in\BZ[t^{\pm1}]$, 
such that 
  \[
\sum_{n\in\BZ}\rank_k(C_n)t^n=p_C(t)\big/\prod_{j=1}^c(1-t^{d_j})\,;
  \]
this is the Hilbert-Serre Theorem, see \cite[4.2]{Sm} or \cite[13.2]{Mat}.  
By \cite[5.5]{Sm}, the order of the pole at $t=1$ of the rational function 
above is equal to the \emph{Krull dimension} of $C$; that is, to the 
supremum of the lengths of chains of homogeneous prime ideals in 
$S'$ containing the annihilator of $C$.   The order of the pole is independent 
of the choice of $S'$, so one gets a well defined notion of Krull dimension 
of $C$ over $S$; let $\dim_SC$ denote this number, and set $\dim S=
\dim_SS$.  By decomposing the rational function above into prime 
fractions, see \cite[\S 2]{Av:msri} for details, one easily obtains 
  \[
\dim_{S}C=\cxy[{}](\rank_k(C_n))\,.
  \]
  \end{chunk}

\begin{theorem}
\label{thm:dghcx}
If $A$ is a DG Hopf $k$-algebra with $S=\Ext{}Akk$
finitely generated as $k$-algebra, and $U$ a DG $A$-module 
with $\hh U$ finite over $k$, then one has
  \[
\injcxy[A]U = \cxy[A]{(U,U)}=\cxy[A]U\leq \cxy[A]k=\dim S\,.
  \]
If, furthermore, $V$ is a DG $A$-module with $\hh V$ finite over $k$, 
then one has
  \[
\dim_S(\Ext{}{A}UV)=\cxy[A]{(U,V)}\le\min\{\cxy[A]{(U,U)},\cxy[A]{(V,V)}\}\,.
  \]
  \end{theorem}

\begin{proof}
The expression for $\cxy[A]{(U,V)}$ comes from Propositions 
\ref{prop:centrality}, \ref{prop:dghfg}, and \ref{ch:HS}.  By Proposition
\ref{prop:dghfg}, $\Ext{}{A}UV$ is a finite module over $\Ext{}{A}UU$ and 
$\Ext{}{A}VV$, whence the upper bounds on $\cxy[A]{(U,V)}$.
For $(k,U)$ and $(U,k)$ they yield
  \[
\injcxy[A]{U}\le\cxy[A]{(U,U)}\ge\cxy[A]{U}\le\cxy[A]{k}\,.
  \]

To prove $\cxy[A](U,U) \leq \cxy[A]U$ we show that the full subcategory 
\[
\mathsf D=\{W\in\dcat A\var \cxy[A](U,W)\leq \cxy[A]U\}
\]
of $\dcat A$ contains $U$.  One evidently has $k\in\mathsf D$, and
\ref{ch:klevel}(1) gives $U\in\thick Ak$, so it suffices to prove that 
$\mathsf D$ is thick.  Closure under direct summands and 
shifts is clear.  An exact triangle $W'\to W\to W''\to \shift W'$ in 
$\dcat A$ yields an exact sequence
  \[
\Ext nAU{W'} \to\Ext nAU{W} \to \Ext nAU{W''}
  \]
for every $n\in\BZ$.  They imply inequalities
  \[
\rank_k\big(\Ext nAU{W}\big)\le
\rank_k\big(\Ext nAU{W'}\big)+
\rank_k\big(\Ext nAU{W''}\big)\,.
  \]
Thus, if $W'$ and $W''$ are in $\mathsf D$, then, Lemma~\ref{lem:artinian} 
gives 
  \[
\cxy[A]{(U,W)}\le\max\{\cxy[A]{(U,W')},\cxy[A]{(U,W'')}\}\le\cxy[A]U\,.
  \]

A similar argument, now using $\Ext{}{A}-U$, yields $\injcxy[A]U=\cxy[A](U,U)$.

The equality $\cxy[A]k =\dim S$ comes from the Hilbert-Serre theorem, see 
\ref{ch:HS}. 
  \end{proof}

To compare $\cxy[A](U,V)$ and $\cxy[A](V,U)$ we need one more lemma.

\begin{lemma}
\label{lem:symmetry}
Let $A$ be a DG Hopf $k$-algebra.  For all $U,V\in\dcat A$ the assignment
$(u\otimes v)\mapsto[\alpha\mapsto(-1)^{|v||\alpha|}u\alpha(v)]$ defines a
morphism
\[
U\otimes_{k}V \to \Hom k{V^{*}}U
\] 
in $\dcat A$; it is an isomorphism when $\rank_{k}\hh V$ is finite.
\end{lemma}

\begin{proof}
It is a routine calculation to verify that the map in question is compatible with 
the DG $A$-module structures. It thus defines for each $X\in\dcat A$ a natural 
morphism $\eta_X\col U\otimes_{k}X\to\Hom k{X^{*}}U$. It is easy to see
that those $X$ for which $\eta_X$ is an isomorphism form a thick subcategory
of $\dcat A$. It contains $k$ and hence contains every DG $A$-module $V$ 
with $\rank_{k}\hh V$ finite; see \ref{ch:klevel}(1).
\end{proof}

The next result extends the equality $\cxy[A]U=\injcxy[A]U$ from
Theorem~\ref{thm:dghcx}.

\begin{theorem}
\label{thm:symmetry}
If $A$ is a DG Hopf $k$-algebra with $\Ext{}Akk$ finitely generated 
as $k$-algebra, and $U$ and $V$ are DG $A$-modules with 
$\hh U$ and $\hh V$ finite over $k$, then
\[
\cxy[A]{(U,V)}=\cxy[A](V,U)\,.
\]
  \end{theorem}

\begin{proof}
The desired assertion results from the following chain of equalities:
\begin{align*}
\cxy[A]{(U,V)} 
			&=\cxy[A](U,V^{**}) \\
			&=\cxy[A](U\otimes_{k}V^{*},k) \\
			&=\cxy[A](k,U\otimes_{k}V^{*}) \\
			&=\cxy[A](k,\Hom kVU) \\
			&=\cxy[A](V,U)\,.
\end{align*}
The first and fourth ones come, respectively, from isomorphisms
$V\simeq V^{**}$ and $U\otimes_{k}V^{*}\simeq\Hom kVU$ 
given by Lemma~\ref{lem:symmetry}. The second and fifth ones
one follow from the adjunction isomorphism \eqref{iso:adjunction}. 
Theorem~\ref{thm:dghcx} supplies the middle link.
  \end{proof}

We  have not yet provided any non-trivial example of DG Hopf 
$k$-algebra $A$ with finitely generated cohomology algebra.  To 
state a general result in this  direction, recall that $A$ is said to
be \emph{cocommutative} if its comultiplication satisfies the 
equality $\Delta=\tau\Delta$, where $\tau\col A\otimes_kA\to
A\otimes_kA$ is defined by $\tau(a\otimes b)=(-1)^{|a||b|}b\otimes a$.

The following result is Wilkerson's main theorem in \cite{Wi}:

\begin{chunk}
  \label{ch:cocommutative}
If $A$ is a cocommutative DG Hopf $k$-algebra with zero differential
and with $\rank_kA$ finite, then the graded $k$-algebra $S=\Ext{}Akk$ is
finitely generated.
  \end{chunk}

In positive characteristic the proof of the theorem depends on the
action of the Steenrod algebra on $S$; the existence of such an
action is another non-trivial result, see \cite{Li}.  The situation is completely
different in characteristic zero, where a celebrated theorem of Hopf
shows that every non-negative cocommutative graded Hopf algebra 
is isomorphic, as an algebra, to the exterior algebra on a vector space 
of finite rank; in this case the cohomology is well known, and is computed next.

\section{Exterior algebras I.}
\label{Exterior algebras. I}

In this section $k$ is a field and $\Lambda$ a DG algebra with $\dd^\Lambda=0$
and with $\Lambda^\natural$ an exterior $k$-algebra on alternating
indeterminates $\xi_{1},\dots,\xi_{c}$ of positive odd degrees.

\begin{proposition}
\label{prop:ext_cohomology}
The graded $k$-algebra $S=\Ext{*}{\Lambda}kk$ is a polynomial
ring over $k$ on commuting indeterminates $\chi_{1},\dots,\chi_{c}$
with $|\chi_{i}|=-(|\xi_{i}|+1)$.
  \end{proposition}

\begin{proof}
One has $\Lambda\cong
k\langle\xi_1\rangle\otimes_{k}{B}$, where $k\langle\xi_j\rangle$
denotes the exterior algebra on the graded vector space $k\xi_i$ and
$B=k\langle\xi_2\rangle\otimes_{k}\cdots\otimes_{k}k\langle\xi_c\rangle$.
By Lemma \ref{lem:Kunneth} and induction, it suffices to treat the case
$\Lambda=k\langle\xi\rangle$.  We use the notation from Construction
\ref{con:composition}.

Set $d=|\xi|$, and let $F$ be a DG $\Lambda$-module, whose underlying
graded $\Lambda$-module has a basis $\{e_{n}\}_{n\ges 0}$ with $|e_{n}| =
(d+1)n$, and with differential defined by the formulas $\dd(e_0)=0$ and
$\dd(e_{n})=\xi e_{n-1}$ for $n\ge1$.  An elementary verification shows
that the $\Lambda$-linear map $\eps\col F\to k$ with $\eps(e_{0})=1$
and $\eps(e_{n})=0$ for $n\ge 1$ is a quasi-isomorphism of DG
$\Lambda$-modules, and thus a semifree resolution of $k$.   It yields
\[
\Ext n{\Lambda}kk = \CH n{\Hom{\Lambda}Fk}=\Hom{\Lambda}Fk^n=
  \begin{cases}
 k\chi^{(i)}&\text{ for } n=(d+1)i\ge0\,;
   \\
 0&\text{ otherwise}\,,
   \end{cases}
\]
where $\chi^{(i)}\col F\to k$ is the chain map of DG $\Lambda$-modules 
defined by $\chi^{(i)}(e_{n})=0$ for $n\ne i$ and $\chi^{(i)}(e_{i})=1$.  
Setting $\wt\chi^{(i)}(e_{n})=e_{n-i}$ for $n\ge i$ and $\wt\chi^{(i)}(e_{n})=0$ 
for $n<i$ one obtains a chain map $\wt\chi^{(i)}\col F\to F$, such that $\eps\wt\chi^{(i)}
=\chi^{(i)}$.  The definition of  composition products yields equalities
  \[
[\chi^{(j)}][\chi^{(i)}]=[\chi^{(j)}\circ\wt\chi^{(i)}]=[\chi^{(i+j)}]\,.
  \]
for all $i,j\ge0$. They show that the isomorphism of graded $k$-vector spaces 
$k[\chi]\to S$, which sends $\chi^i$ to $[\chi^{(i)}]$ for each $i\ge0$, is a 
homomorphism of graded $k$-algebras.
 \end{proof}

\begin{remark}
\label{rem:ext_Hopf}
The universal property of exterior algebras guarantees that there is a 
unique homomorphism of graded $k$-algebras 
$\Delta\col \Lambda\to\Lambda\otimes_{k}\Lambda$ with
\[
\Delta(\xi_{i}) = \xi_{i}\otimes 1 + 1\otimes \xi_{i} \quad\text{for $1\leq i\leq c$}.
\]
It evidently satisfies $\Delta(\eps^A\otimes A)=\idmap^A=
\Delta(A\otimes_k\eps^A)$, and so is the comultiplication of a graded
Hopf algebra structure on $\Lambda$.  
  \end{remark}

In view of the preceding proposition and remark, Theorems \ref{thm:dghcx} and \ref{thm:symmetry} yield:

\begin{theorem}
\label{thm:ext_cx}
\pushQED{\qed}
For DG $\Lambda$-modules $U,V$ with $\hh U$ and 
$\hh V$ finite over $k$ one has
\[
\cxy[A]{(V,U)}=\cxy[A]{(U,V)}\le
\cxy[\Lambda]{(U,U)}=\injcxy[\Lambda]U=\cxy[\Lambda]U\leq \cxy[\Lambda]k=c\,.
  \qedhere
\]
\end{theorem}

\section{Complete intersection local rings I.}
\label{Complete intersection local rings.I}

In this section $(R,\fm,k)$ denotes a local ring.  When $R$ is complete 
intersection (the definition is recalled below) we open a path to an exterior 
algebra, and use it to transport results from Section \ref{Exterior algebras. I}.

The \emph{embedding dimension} of $(R,\fm,k)$ is the number $\edim R
=\rank_{k}(\fm/\fm^{2})$; by Nakayama's Lemma, it is equal to the minimal
number of generators of $\fm$.  

\begin{construction}
\label{con:koszul_k}
Choose a minimal set of generators $r_1,\dots,r_e$ of $\fm$, and let 
$K$ denote the Koszul complex on $r_1,\dots,r_e$.
The functor $-\otimes_RK$ preserves quasi-isomorphisms
of complexes of $R$-modules, so it defines an exact functor 
  \[
\SK\col\dcat R\to\dcat K\,.
  \]
  \end{construction}

  \begin{lemma}
   \label{lem:koszul_k}
For each homologically finite complex of $R$-modules $M$ one has:
  \begin{align*}
\cxy M=\cxy[K]{\SK(M)}
  \quad\text{and}\quad
\injcxy M=\injcxy[K]{\SK(M)}\,.
  \end{align*}
  \end{lemma}

\begin{proof}
Let $G\to M$ be a semifree resolution over $R$, and note that the induced
morphism $G\otimes_RK\to M\otimes_RK$ then is a semifree resolution over $K$.
The expression for $\cxy M$ now follows from the isomorphisms of $k$-vector 
spaces
\begin{align*}
\Ext n{K}{M\otimes_RK}{k}
	&=\CH n{\Hom K{G\otimes_RK}{k}} \\
	&\cong\CH n{\Hom R{G}{k}} \\
	&=\Ext n{R}{M}{k}\,.
\end{align*}

For the second equality, choose a semifree resolution $F\to k$
over $K$.  One has an isomorphism $\SK M\simeq
\Hom RK{\shift^{e}M}$ of DG $K$-modules, where $e=\edim R$.
For each $n\in\BZ$ adjunction gives the isomorphism of $k$-vector 
spaces below, and the last equality holds because $F$ is semifree 
over $R$:
  \begin{align*}
\Ext n Kk{M\otimes_RK}
   &=\CH n{\Hom K{F}{\Hom RK{\shift^{e}M}}} \\
   &\cong\CH n{\Hom R{F}{\shift^{e}M}} \\
   &=\Ext {n-e}{R}{k}{M}\,.
     \qedhere
  \end{align*}
  \end{proof}

Krull's Principal Ideal Theorem implies
$\edim R\geq \dim R$, where $\dim R$ denotes the Krull dimension of $R$.
The \emph{codimension} of $R$ is the number $\codim R=\edim R-\dim R$.
Rings of codimension zero are called \emph{regular}.  Cohen's  Structure
Theorem shows that the $\fm$-adic completion $\wh R$ of $R$ admits a
surjective homomorphism of rings $\varkappa\col (Q,\fq,k)\to \wh R$, with $Q$
a regular local ring; see \cite[\S29]{Mat}.  Such a homomorphism is
called a \emph{Cohen presentation} of $\wh R$.  One has $\dim Q\geq
\edim R$, and equality is equivalent to $\Ker(\varkappa)\subseteq\fq^{2}$;
such a presentation is said to  be \emph{minimal}.

Any Cohen presentation $\varkappa$ can be refined to a minimal one.  Indeed,
choose elements $q_{1},\dots,q_{n}$ in $\fq$ mapping to a $k$-basis of
$(\Ker(\varkappa) + \fq^{2})/\fq^{2}$.  One then has
\[
\codim Q/(\bsq) = (\edim Q - n)-(\dim Q - n) =\codim Q=0\,,
\]
so the local ring $Q/(\bsq)$ is regular along with $Q$.  The map $\varkappa$
factors through a homomorphism $Q/(\bsq)\to\wh R$, and the latter is
a minimal Cohen presentation.

\begin{construction}
\label{con:koszul_j}
Let $\iota\col R\to\wh R$ be the completion map and choose a
minimal Cohen presentation $\varkappa\col(Q,\fq,k)\to\wh R$.

Choose a minimal set of generators $r_1,\dots,r_e$ of $\fm$ and then pick
in $Q$ elements $q_1,\dots,q_e$ with $\varkappa(q_i)=\iota(r_i)$ for
$i=1,\dots,e$.  Let $K$ denote the Koszul complex on $r_1,\dots,r_e$,
and $E$ be the one on $q_1,\dots,q_e$.  The definition of Koszul
complexes allows one to identify the DG algebras ${\wh R}\otimes_{R} K$
and ${\wh R}\otimes_{Q}E$.

Choose a minimal set of generators $\{f_1,\dots,f_c\}$ of the ideal
$\Ker\varkappa$, let $A$ be the Koszul complex on this set, and let
$\pi\col A\to \wh R$ denote the canonical projection.

Set $\Lambda=A\otimes_{Q} k$, and note that this is a DG algebra with
zero differential, and its underlying graded algebra is the exterior
algebra $\textstyle{\bigwedge}_k(A_1\otimes_Qk)$.
  \end{construction}

The ring $R$ is said to be \emph{complete intersection} if in some Cohen 
presentation $\varkappa\col Q\to \wh R$ the ideal $\Ker\varkappa$ can be
generated by a $Q$-regular sequence.  When this is the case, the kernel
of \emph{every} Cohen presentation is generated by a regular sequence,
and for each minimal presentation such a sequence consists of $\codim
R$ elements.

In the next three lemmas
\emph{$(R,\fm,k)$ denotes a complete intersection local ring}.

  \begin{lemma}
   \label{lem:koszul_j}
The following maps are quasi-isomorphisms of DG algebras:
 \[
\xymatrixcolsep{2.5pc}
\xymatrix{
K=R\otimes_RK\ar@{->}[r]^-{\iota\otimes_RK}
        & {\wh R}\otimes_{R} K = {\wh R}\otimes_{Q} E
              \ar@{<-}[r]^-{\pi\otimes_Q E}
	& A\otimes_{Q} E \ar@{->}[r]^-{A\otimes_Q\eps^E}
	& A\otimes_{Q}k=\Lambda\,.
	}
 \]
  \end{lemma}

\begin{proof}
As $\wh R$ is flat over $R$, one can identify $\hh{\iota\otimes_RK}$ 
with the map
\[
\iota\otimes_R\hh K\col R\otimes_R\hh K\to\wh R\otimes_R\hh K\,,
\]
which is bijective because $\hh K$ is a direct sum of shifts of $k$.

The sequence $f_1,\dots,f_c$ is $Q$-regular because $R$ is complete
intersection, hence $\pi$ is a quasi-isomorphism, and then so is
$\pi\otimes_QE$.

Since the ring $Q$ is regular and the elements $q_1,\dots,q_e$ minimally
generate $\fq$, they form a $Q$-regular sequence, so  $\eps^{E}$ 
is a quasi-isomorphism, and then so is $A\otimes_Q\eps^E$.
\end{proof}

In view of \ref{chu:quism}, the quasi-iso\-morph\-isms in Lemma 
\ref{lem:koszul_j} define an equivalence of categories
$\SJ\col\dcat K\to\dcat\Lambda$.  Lemmas \ref{lem:koszul_k} and 
\ref{lem:quism} 
then give:

  \begin{lemma}
   \label{lem:koszul_cx}
    \pushQED{\qed}
For every homologically finite complex of $R$-modules $M$ one has
  \[
\cxy M=\cxy[\Lambda]{\SJ\SK(M)}
  \quad\text{and}\quad
\injcxy M=\injcxy[\Lambda]{\SJ\SK(M)} \,.
   \qedhere
  \]
  \end{lemma}
  
  \begin{lemma}
   \label{lem:koszul_jk}
In $\dcat\Lambda$ there is an isomorphism $\SJ\SK(k) \simeq
\bigoplus_{i\ges 0} \shift^{i}k^{\binom ci}$.
  \end{lemma}

  \begin{proof}
Choose bases $x_1,\dots,x_c$ of $A_1$ over $Q$ with $\dd(x_i)=f_i$
for $1\leq i\leq c$ and $y_1,\dots,y_e$ of $E_1$ over $Q$ with 
$\dd(y_j)=q_j$ for $1\leq j\leq e$; Thus $f_{i}=\sum_{j=1}^{d}b_{ij}q_{j}$ 
with $b_{ij}\in Q$.  There is a unique homomorphism $\alpha\col A\to E$ 
of DG algebras over $Q$ with
  \[
\alpha(x_i)=\sum_{j=1}^{e}b_{ij}y_j
  \quad\text{for}\quad 1\leq i\leq c\,.
  \]
It appears in a commutative diagram of DG algebras 
 \[
\xymatrixcolsep{.5pc}
\xymatrixrowsep{2.5pc}
\xymatrix{
K\ar@{=}[r]
        & R\otimes_RK\ar@{->}[rrr]^-{\iota\otimes_RK}
          \ar@{->}[d]^-{\eps^R\otimes_RK}
        &&& {\wh R}\otimes_{R} K \ar@{=}[r]
          \ar@{->}[d]^-{\eps^{\wh R}\otimes_QE}
        & {\wh R}\otimes_{Q} E \ar@{<-}[rrr]^-{\pi\otimes_Q E}
          \ar@{->}[d]^-{\eps^{\wh R}\otimes_QE}
	&&& A\otimes_{Q} E \ar@{->}[rrr]^-{A\otimes_Q\eps^E}
	  \ar@{->}[d]^-{\alpha\otimes_QE}
	&&& A\otimes_{Q} k\ar@{=}[r]
	  \ar@{->}[d]^-{\alpha\otimes_Qk}
& \Lambda
  \\
\SK(k)\ar@{=}[r]
& k\otimes_RK\ar@{=}[rrr]
&&& k\otimes_RK\ar@{=}[r]
&   k\otimes_QE\ar@{<-}[rrr]^-{\eps^E\otimes_Q E}_-{\simeq}
&&& E\otimes_QE\ar@{->}[rrr]^-{E\otimes_Q\eps^E}_-{\simeq}
&&& E\otimes_Qk
}
 \]
where $\eps^E\otimes_Q E$ and $E\otimes_Q\eps^E$ are 
quasi-isomorphisms because $\eps^E$ is one.

The map $\epsilon^E\otimes_RK$ defines the action of $K$ on $\SK(k)$,
and $\alpha\otimes_Qk$ that of $\Lambda$ on $E\otimes_Qk$.  The 
commutativity of the diagram implies and $\SJ\SK(k)\simeq E\otimes_Qk$ 
in $\dcat\Lambda$, see \ref{chu:quism}.  The minimality 
of the Cohen presentation $\varkappa$ means that each $b_{ij}$ is in $\fq$, 
so one has $\SK(k)\cong\bigoplus_{i\ges 0}\shift^{i}k^{\binom ci}$
as DG $K$-modules, and hence in $\dcat\Lambda$ one gets
  \[
\SJ\SK(k)
\simeq\SJ\bigg(\bigoplus_{i\ges 0}\shift^{i}k^{\binom ci}\bigg)
\simeq\bigoplus_{i\ges 0}\shift^{i}\SJ(k)^{\binom ci}
\simeq\bigoplus_{i\ges 0}\shift^{i}k^{\binom ci}\,,
 \]
where the last isomorphism holds because one has $\SJ(k)\simeq k$ by
Lemma \ref{lem:quism}.
  \end{proof}

We come to the main result of this section.  All its assertions are 
known:  the first equality is proved in \cite[5.3]{Av:vpd}, the inequality 
follows from \cite[4.2]{Gu},
and the second equality from \cite[Thm.~6]{Ta}.   The point here is that
they are deduced, in a uniform way, from the corresponding relations
for DG modules over exterior algebras, established in Theorem 
\ref{thm:ext_cx}, which ultimately are much simpler to prove.  

  \begin{theorem}
   \label{thm:ci_cx}
If $(R,\fm,k)$ is complete intersection, then for every complex of $R$-modules
$M$ with $\hh M$ finitely generated the following inequalities hold:
 \[
\injcxy M=\cxy M\leq \cxy k=\codim R\,.
  \]
    \end{theorem}

\begin{proof}
The isomorphism of Lemma \ref{lem:koszul_jk} implies 
$\cxy[\Lambda]k=\cxy[\Lambda]{\SJ\SK(k)}$.
Now Lemma \ref{lem:koszul_cx} translates part of Theorem
\ref{thm:ext_cx} into the desired statement.
  \end{proof}

The equalities in the theorem may fail when $R$ is not complete 
intersection:

\begin{remark}
\label{ex:js}
Gulliksen \cite[2.3]{Gu1} proved that the condition $\cxy  k<\infty$ 
characterizes  local complete  intersection rings among all local rings.   

Jorgensen and \c{S}ega \cite[1.2]{JS} construct
Gorenstein local $k$-algebras $R$ with $\rank_kR$ finite and modules 
$M$ with $\{\cxy M,\injcxy M\}=\{1,\infty\}$, in one order or the other.
\end{remark}

\begin{remark}
The unused portion of Theorem \ref{thm:ext_cx} suggests that
the relations
 \[
\cxy{(N,M)}=\cxy{(M,N)}\le\cxy{(M,M)}=\cxy M
  \]
might hold also over complete intersections.  They do, see \cite[5.7]{AB}, 
but we know of no simple way to deduce them from Theorem \ref{thm:ext_cx}. 
 \end{remark}

\section{Projective levels}
\label{Projective levels}

In this section $A$ is a DG algebra.  We collect some results 
on $A$-levels of DG modules, which are reminiscent of theorems on 
projective dimension for modules.

  \begin{chunk}
\label{ch:Alevel}
For any DG $A$-module $U$ the following hold.
\begin{enumerate}[{\rm(1)}]
\item
One has $\level AAU\leq l+1$ if and only if $U$ is isomorphic in $\dcat A$ to
a direct summand of some semifree DG module $F$ with a semifree filtration 
having $F^{l}=F$ and every $F^n/F^{n-1}$ of finite rank over $A^\natural$;
see \ref{ch:Ext}.
  \item
When $A$ is non-negative, $\HH0A$ is a field, and $\HH iU=0$ for all 
$i\ll0$, one has $\level AAU<\infty$ if and only if $\Tor{}AkU$ is finite over $k$.
  \item
When $A$ has zero differential and is noetherian, and the graded 
$A$-module $\hh U$ is finitely generated, the following inequalities hold:
    \[
\level AAU\le\pd_A\hh U+1\le\gldim A+1\,.
    \]
  \end{enumerate}
Indeed, (1) is proved in \cite[4.2]{ABIM}, (2) in \cite[4.8]{ABIM}, and 
(3) in \cite[5.5]{ABIM}.
    \end{chunk}

The preceding results pertain to homological algebra, in the sense that 
they do not depend on the structure of the ring $A^\flat=\bigoplus_{n\in\BZ} 
A_n$.  For algebras over fields the next result contains as a special case the 
New Intersection Theorem, see~\cite{Rb}, a central result in commutative 
algebra, and effectively belongs to the latter subject.

\begin{chunk}
\label{ch:nit}
Let $A$ be a DG algebra with zero differential and $U$ a DG $A$-module.

If the ring $A^\flat$ is commutative, noetherian, and is an algebra over a field, then
 \[
\level AAU\geq \height I+1\,,
 \] 
where $I$ is the annihilator of the $A^\flat$-module 
$\bigoplus_{n\in\BZ}\HH nU$; see \cite[5.1]{ABIM}.
\end{chunk}

\section{Exterior algebras II.}
\label{Exterior algebras.II}

In this section $k$ is a field, and $\Lambda$ a DG algebra with 
$\dd^\Lambda=0$ and $\Lambda^\natural$ an exterior $k$-algebra
on alternating indeterminates of positive odd degrees $d_1,\dots,d_c$.

We recall a version of the Bernstein-Gelfand-Gelfand equivalence
from~\cite[7.4]{ABIM}: 

\begin{chunk}
\label{ch:bgg}
Let $S$ be a DG algebra with $\dd^S=0$ and $S^\natural$ a polynomial 
ring over $k$ on commuting indeterminates of degrees 
$-(d_1+1),\dots,-(d_c+1)$; set $d=d_1+\cdots+d_c$.

There exist exact functors 
$\dcat{\Lambda}\xra{\SH}\dcat{S}\xra{\ST}\dcat{\Lambda}$ 
inducing inverse equivalences 
  \begin{equation}
    \label{eq:bgg_equiv}
 \begin{gathered}
\xymatrixcolsep{2.5pc} 
\xymatrixrowsep{1.5pc} 
\xymatrix{
\thick{\Lambda}{k}\ar@{->}[r]<1.2ex>^-{\SH}_-{\equiv}
\ar@{}[d]|-{\vgu}
&\thick{S}{S}\ar@{->}[l]<1ex>^-{\ST}
\ar@{}[d]|-{\vgu}
\\
\thick{\Lambda}{\Lambda}\ar@{->}[r]<1.2ex>^-{\SH}_-{\equiv}
&\thick{S}{k}\ar@{->}[l]<1ex>^-{\ST}
}
 \end{gathered}
\end{equation}
such that in $\dcat S$ and $\dcat\Lambda$, respectively, there are 
isomorphisms
    \begin{equation}
      \label{eq:bgg_values}
\begin{aligned}
\SH(\Lambda)&\simeq \shift^d k
\\
\SH(k)&\simeq S
  \end{aligned}
\qquad\text{and}\qquad
\begin{aligned}
\ST(k)&\simeq\shift^{-d}\Lambda
\\
\ST(S)&\simeq k\,.
  \end{aligned}
    \end{equation}
  \end{chunk}

Let $\dcatf{\Lambda}$ be the full subcategory of $\dcat{\Lambda}$
whose objects are the DG modules $U$ with $\hh U$ finite
over $\Lambda$, and $\dcatf S$ the corresponding
subcategory of $\dcat S$.  The next proposition 
refines the equivalences in \eqref{eq:bgg_equiv}; see
Remark \ref{rem:refinement}.

\begin{proposition}
\label{lem:bgg}
For each $U\in\dcatf{\Lambda}$ and each $M\in\dcatf S$ one has
  \[
\cxy[\Lambda]U = \dim_{S}\hh{\SH U}
  \quad\text{and}\quad
\dim_{S}M=\cxy[\Lambda]{\hh{\ST M}}\,,
  \]
and for each $n\in\BN$ the functors $\SH$ and $\ST$ induce inverse 
equivalences
  \begin{equation}
    \label{eq:bgg_equiv1}
 \begin{gathered}
\xymatrixcolsep{2.5pc} 
\xymatrixrowsep{1.5pc} 
\xymatrix{
\dcatf{\Lambda}\ar@{->}[r]<1.2ex>^-{\SH}_-{\equiv}
\ar@{}[d]|-{\vgu}
&\dcatf{S}\qquad\ar@{->}[l]<1ex>^-{\ST}
\ar@{}[d]|-{\vgu\qquad}
\\
\{U\in \dcatf{\Lambda} \mid \cxy[\Lambda]U \leq n\}
\ar@{->}[r]<1.2ex>^-{\SH}_-{\equiv}
&\{M\in \dcatf{S}\mid \dim_{S}\hh M \leq n \}
\ar@{->}[l]<1ex>^-{\ST}\,.
}
 \end{gathered}
\end{equation}
    \end{proposition}

\begin{proof}
One has $\dcatf\Lambda=\thick{\Lambda}k$ by \ref{ch:klevel}(1).  As the
polynomial ring $S$ has finite global dimension, \ref{ch:Alevel}(3) gives
$\dcatf S=\thick SS$.  Thus, the top row of \eqref{eq:bgg_equiv} gives
inverse equivalences between $\dcatf\Lambda$ and $\dcatf S$.  It is easy
to verify that the subcategories in the lower row are thick, so it suffices
to compute $\cxy[\Lambda] U$ and $\dim_{S}\hh M$.

The equivalence $\SH$ gives the first isomorphism in the chain
\[
\Rhom{\Lambda}kU\simeq\Rhom{S}{\SH(k)}{\SH(U)}
\simeq\Rhom{S}{S}{\SH(U)}\simeq\SH(U)\,;
\]
the second one comes from  \eqref{eq:bgg_values}, the third is clear.
In homology, one obtains 
\[
\Ext n{\Lambda}kU\cong\CH n{\SH(U)}
\]
for every $n\in\BZ$, whence the second equality in the following 
sequence; the first one comes from Theorem \ref{thm:ext_cx}, 
the third from the Hilbert-Serre Theorem, see~\ref{ch:HS}:
\[
\cxy[\Lambda] U = 
\injcxy[\Lambda] U = 
\cxy[{}](\rank_{k}\HH n{\SH(U)})_{n\in\BZ}
= \dim_{S} \hh{\SH(U)}\,.
\]

Finally, from the already proved assertions one obtains
  \[
\dim_{S}\hh M=\dim_{S}{\hh{\SH\ST(M)}}=\cxy[\Lambda]{\ST(M)}\,.
 \qedhere
  \]
  \end{proof}

  \begin{remark}
    \label{rem:refinement}
For $U\in\dcatf\Lambda$ the definition of complexity shows that 
$\cxy[\Lambda]U\leq 0$ is equivalent to the finiteness of the 
$k$-vector space $\Ext{}{\Lambda}Uk$.  It is the graded $k$-dual 
of $\Tor{}{\Lambda}Uk$, so by \ref{ch:Alevel}(2) its finiteness 
is equivalent to $U\in\thick{\Lambda}{\Lambda}$.  On the other
hand, $M\in\dcatf S$ has $\dim_S\hh M\le0$ if and only if $M$ is in
$\thick Sk$; see \ref{ch:klevel}(1).  Thus, for $n=0$ diagram 
\eqref{eq:bgg_equiv1} reduces to \eqref{eq:bgg_equiv}.
  \end{remark}

\begin{theorem}
\label{thm:ext_main}
If $U$ is a DG $\Lambda$-module with $0<\rank_{k}\hh U<\infty$, then
one has
\[
\card\{n \in\BZ\mid \HH nU\ne 0\}\ge
\level {\Lambda}kU \geq c - \cxy[\Lambda] U +1\,.
\]
\end{theorem}

\begin{proof}
The inequality on the left comes from \ref{ch:klevel}(3). The one on the right
results from the following chain of (in)equalities
  \begin{align*}
\level {\Lambda}kU 
	&= \level{S}{\SH(k)}{\SH(U)} \\
	&=\level{S}{S}{\SH(U)}\\
	&\geq \height({\Ann_{S}\hh{\SH(U)}}) +1 \\
	&=\dim S - \dim_{S}\hh{\SH(U)} +1 \\
	&=\dim S - \cxy[\Lambda]U + 1  \\
        &=c- \cxy[\Lambda]V + 1\,.
  \end{align*}
The first one holds because $\SH$ is an equivalence, see \ref{ch:levels}(2), the second comes from \eqref{eq:bgg_values}, and the third from \ref{ch:nit}; the remaining equalities hold by \cite[Exercise 5.1]{Mat}\footnote{solved on page 288}, by Lemma~\ref{lem:bgg}, and because $\dim S=c$ holds.
  \end{proof}

The inequalities in Theorem~\ref{thm:ext_main} are tight.

\begin{example}
\label{ex:ext_main}
For each integer $i$ with $0\leq i\leq c$ and $\Lambda^{(i)}=
\Lambda/(\xi_{1},\dots,\xi_{i})$ one has
\[
\card\{n \in\BZ\mid \HH n{\Lambda^{(i)}}\ne 0\}= c-i+1\, 
\quad\text{and}\quad \cxy \Lambda^{(i)}= i\,.
\] 

Indeed, the first one is clear.  For $\Gamma^{(i)}=\Lambda/(\xi_{i+1},\dots,\xi_{n})$ 
one has
  \[
\Ext n{\Lambda}{\Lambda^{(i)}}k
\cong\Ext n{\Gamma^{(i)}\otimes_k{\Lambda^{(i)}}}{\Lambda^{(i)}}k
\cong\Ext n{\Gamma^{(i)}}kk\,,
  \]
because $\Lambda$ is isomorphic to $\Gamma^{(i)}\otimes_k\Lambda^{(i)}$ 
as $k$-algebras.  Theorem \ref{thm:ext_cx} gives $\cxy[\Gamma^{(i)}]k=i$.
  \end{example}

\section{Complete intersection local rings II.}
\label{Complete intersection local rings.II}

We finish the paper with a new proof of \cite[11.3]{ABIM}.  The original 
argument uses a reduction, constructed in \cite{Av:vpd}, to a complete 
intersection ring of smaller codimension and a bounded complex of 
free modules over it, then refers to the main theorem of \cite{ABIM}, 
which applies to such complexes over arbitrary local rings.   Here we 
describe a direct reduction to results on exterior algebras, established 
in Section \ref{Exterior algebras.II}.

\begin{theorem}
\label{thm:ci_main}
If $(R,\fm,k)$ is a complete intersection local ring and $M$ a complex
of $R$-modules with $\hh M$ finite and nonzero, then one has inequalities
\[
\sum_{n\in\BZ}\lol R{\HH nM}\ge\level RkM \geq \codim R - \cxy M  +1\,.
\]
\end{theorem}

\begin{proof}
The inequality on the left is a special case of \ref{ch:klevel}(2).

Set $c=\codim R$.  Let $K$ be the Koszul complex on a minimal set of
generators for $\fm$ and $\Lambda$ the DG algebra with zero differential,
and with underlying graded algebra an exterior algebra over $k$ on $c$
generators of degree one.  Construction \ref{con:koszul_k} and Lemma
\ref{lem:koszul_j} provide exact functors of triangulated categories
    \begin{equation*}
    \xymatrix{
\dcat R \ar@{->}[r]^{\SK}
     &\dcat K\ar@{->}[r]^-{\SJ}_-\equiv
     &\dcat \Lambda
     }
  \end{equation*}
the second of which is an equivalence.  One then has the following relations
\begin{align*}
\level RkM 
           &\ge \level {\Lambda}{\SJ\SK(k)}{\SJ\SK(M)} \\
           &= \level {\Lambda}k{\SJ\SK(M)} \\
           &\geq c - \cxy[\Lambda](\SJ\SK(M))+1 \\
           &= \codim R - \cxy M+1
\end{align*}
where the inequalities are given by \ref{ch:levels}(2) and by Theorem \ref{thm:ext_main},
while the equalities come from \ref{ch:levels}(1) and Lemmas \ref{lem:koszul_jk},
and from Lemma \ref{lem:koszul_cx}.
  \end{proof}

As in Theorem \ref{thm:ext_main} the inequalities in Theorem~\ref{thm:ci_main} 
are tight.

\begin{example}
\label{ex:cimain}
Let $k$ be a field and set $R=k[x_{1},\dots,x_{c}]/(x_{1}^{2},\dots,x_{c}^{2})$.
For every integer $i$ with $0\leq i\leq c$, set $R^{(i)}=R/(x_{1},\dots,x_{i})$.
One then has equalities
\[
\lol R{R_{i}}= c-i+1\,,\quad\codim R=c\,, \quad\text{and}\quad \cxy R^{(i)} = i\,.
\] 

Indeed, the first two are clear.  For $Q^{(i)}=R/(x_{i+1},\dots,x_{n})$ one has
  \[
\Ext n{R}{R^{(i)}}k
\cong\Ext n{Q^{(i)}\otimes_k{R^{(i)}}}{R^{(i)}}k
\cong\Ext n{Q^{(i)}}kk\,,
  \]
because $R$ is isomorphic to $Q^{(i)}\otimes_k R^{(i)}$ as $k$-algebras.
Theorem \ref{thm:ci_cx} gives $\cxy[Q^{(i)}]k=i$.
\end{example}

The number on the right hand side of the formula in Theorem \ref{thm:ci_main} 
is finite whenever the module $M$ has finite complexity.  However, this condition 
alone does not imply the inequalities in the theorem, even when
the ring $R$ is Gorenstein.

\begin{example}
\label{ex:sqzero}
Let $k$ be a field and $c$ an integer with $c\ge3$.   The ring 
  \[
R=\frac{k[x_{1},\dots,x_{c}]}
{\sum_{h=1}^{c-1}(x_h^2-x_{h+1}^2)+\sum_{1\les i< j\les c}(x_ix_j)}
  \]
is Gorenstein, but not complete intersection, and its module
$M=R$ has
\[
\lol RM = \level RkM = 3 < c-0+1=\codim R-\cxy M+1\,.
\]
\end{example}

It is worth noting that a very special case of Theorem \ref{thm:ci_main} 
was initially discovered when studying actions of
finite elementary abelian groups on finite CW complexes:

  \begin{remark}
Let $k$ be a field of positive characteristic $p$ and $G$ an elementary 
abelian $p$-group of rank $c$. The group algebra $kG$ is isomorphic to
$k[x_{1},\dots,x_{c}]/(x_{1}^{p},\dots,x_{c}^{p})$, which is a complete
intersection of codimension $c$. Over $kG$ Theorem \ref{thm:ci_main}
was proved by Carlsson \cite{Ca:inv} for $M^\natural$ of finite rank over 
$A^\natural$ and $p$ equal to $2$, and for general $M$ and $p$ by 
Allday, Baumgartner, and Puppe; see~\cite[4.6.42]{AP}.
  \end{remark}


\begin{thebibliography}{99}

\bibitem{AP}
C.~Allday, V.~Puppe,
\textit{Cohomological methods in transformation groups},
Cambridge Stud. Adv. Math. \textbf{32},
Cambridge Univ. Press, Cambridge, 1993.

\bibitem{Av:vpd}
L.~L.~Avramov,
\textit{Modules of finite virtual projective dimension},
 Invent. Math. \textbf{96} (1989), 71--101.

\bibitem{Av:msri}
L.~L.~Avramov,
\textit{Homological asymptotics of modules over local rings},
Commutative algebra (Berkeley, CA, 1987),  
Math. Sci. Res. Inst. Publ., \textbf{15}, Springer, New York, 1989;
33--62.

\bibitem{AB}
L.~L.~Avramov, R.-O.~Buchweitz, 
\textit{Support varieties and cohomology over complete intersections},
  Invent. Math.  \textbf{142}  (2000),  285--318.

\bibitem{ABIM}
L.~L.~Avramov, R.-O.~Buchweitz, S.~B.~Iyengar, C.~Miller,
\emph{Homology of finite free complexes}, Adv. Math. \textbf{223} (2010), 1731--1781.

\bibitem{AH}
L.~L.~Avramov, S.~Halperin,
\textit{Through the looking glass: a dictionary between rational homotopy 
theory and local algebra},
Algebra, algebraic topology and their interactions (Stockholm, 1983),  
Lecture Notes in Math. \textbf{1183}, Springer, Berlin, 1986;  1--27

\bibitem{Ca:inv}
     G.~Carlsson,
     \textit{On the homology of finite free $(\BZ/2)^k$-complexes},
     Invent. Math. \textbf{74} (1983), 139--147.

\bibitem{Gu}
T. H. Gulliksen,  
\textit{A change of ring theorem with applications to Poincar\'e series and intersection multiplicity}, 
Math. Scand. \textbf{34} (1974), 167--183. 

\bibitem{Gu1}
T. H. Gulliksen,  
\textit{On the deviations of a local ring}, 
Math. Scand. \textbf{47} (1980), 5--20. 

\bibitem{JS}
D.~A.~Jorgensen, L. M. \c Sega, 
\textit{Asymmetric complete resolutions and vanishing of Ext over Gorenstein rings},
Int. Math. Res. Not.  \textbf{2005},  3459--3477.

\bibitem{Li}
A.~Liulevicius, 
\textit{The factorization of cyclic reduced powers by secondary 
cohomology operations}, Mem. Amer. Math. Soc. \textbf{42},
Amer. Math. Soc., Providence, RI, 1962. 

\bibitem{Mac}
S. Maclane,  
\textit{Homology}, Grundlehren math. Wissenschaften, 
\textbf{114}, Springer, Berlin, 1963.

\bibitem{Mat}
H. Matsumura,  
\textit{Commutative ring theory}, Cambridge Stud. Adv. Math. \textbf{8}
Cambridge Univ. Press, Cambridge, 1986.

\bibitem{Rb}
     P.~C.~Roberts,
\textit{Multiplicities and Chern classes in local algebra},
Cambridge Tracts  Math. \textbf{133},
Cambridge Univ. Press, Cambridge, 1998.

\bibitem{Sm}
W.~Smoke, 
\textit{Dimension and multiplicity over graded algebras},
J. Algebra \textbf{21} (1972), 149--173.

\bibitem{Ta}
J. Tate, 
\textit{Homology of noetherian rings and local rings}, 
Illinois J. Math. \textbf{1} (1957), 14--25

\bibitem{Wi}
C.~Wilkerson, 
\textit{The cohomology algebras of finite-dimensional Hopf algebras},
Trans. Amer. Math. Soc. \textbf{264} (1981), 137--150.

  \end{thebibliography}
\end{document}